\documentclass[12pt]{amsart}
\usepackage{amssymb,amsmath}
\usepackage{amsthm}
\usepackage{amsfonts}
\usepackage[dvips]{graphicx}
\usepackage{enumerate}

\newcommand{\Cd}{\mathcal{C}}
\newcommand{\R}{\mathbb{R}}
\newcommand{\Z}{\mathbb{Z}}
\newcommand{\Q}{\mathbb{Q}}
\newcommand{\C}{\mathbb{C}}
\newcommand{\D}{\mathbb{D}}
\newcommand{\Hy}{\mathbb{H}}
\newcommand{\re}{\mathfrak{Re}}
\newcommand{\im}{\mbox{im}\;}
\renewcommand{\Im}{\mathfrak{Im}}
\newcommand{\K}{\mathfrak{K}}
\renewcommand{\epsilon}{\varepsilon}

\newcommand{\minus}{\smallsetminus}

\newtheorem{theorem}{Theorem}[section]

\newtheorem{prob}{Problem}

\newtheorem{lema}[theorem]{Lemma}
\newtheorem{prop}[theorem]{Proposition}

\newtheorem{claim}{Claim}

\title{On the geodesic flow of surfaces of nonpositive curvature}
\author{Federico Rodriguez Hertz}
\address{IMERL, Montevideo, Uruguay}
\email{frhertz@fing.edu.uy}
\date{\today}
\begin{document}

\begin{abstract}
Let $S$ be a surface of nonpositive curvature of genus bigger than 1 (i.e. not the torus). We prove that any flat strip in the surface is in fact a flat cylinder. Moreover we prove that the number of homotopy classes of such flat cylinders is bounded. 
\end{abstract}

\maketitle
\begin{section}{\bf Introduction}
Let $S$ be a surface of genus bigger than one and let $g$ be a metric of nonpositive curvature in $S$. It is known that if the metric has negative curvature then the geodesic flow of the surface is Anosov. Moreover, it is not hard to prove that if no geodesic stays all the time in the region of zero curvature then the geodesic flow is also Anosov. So, we want to deal with the case that there is a geodesic staying in zero curvature. Let us identify the universal covering of $S$ with $\D$, the unit disc in the complex plane. $\D$ with the induced metric by the covering projection, $p:\D\to S$ becomes a Hadamard manifold, i.e. a simply connected manifold with nonpositive curvature. It is not hard to prove that if two geodesics remain at bounded distance in $\D$ then they bound a flat strip $B$, i.e. an open set bounded by the two geodesics where the curvature vanishes, or equivalently, there exists an isometry $h:L\to\D$ where $L=\{z\in\C\; \mbox{s.t.}\; 0\leq\re(z)\leq \alpha\}$ for some $\alpha>0$ sending $\re(z)=0$ to one of the geodesics and $\re(z)=\alpha$ to the other one. In \cite{le}, it is proven that if there is no such flat strip then the geodesic flow is expansive and conjugated with the geodesic flow of a hyperbolic surface. Moreover, it is proven certain kind of stability of the geodesic flow. On the other hand, it is not hard to prove that if the set of orbits in the unit fiber bundle of the surface, associated with the geodesics lying in the zero curvature region, has Liouville measure zero then the geodesic flow is ergodic (it is in fact Bernoulli and has non zero Lyapunov exponents almost everywhere). These facts and the theorem we will prove motivates the following problem:
\begin{prob}
Call $K=\{x\in S\;\mbox{s.t.}\; \K(x)=0\}$ where $\K$ is the curvature of $g$. Let $\gamma$ be a geodesic in $S$. Then if $\gamma\subset K$ then $\gamma$ is a closed geodesic. Moreover there are only finitely many homotopy classes of such geodesics.
\end{prob} 
In this paper we shall solve the problem for geodesics associated with flat strips.
\begin{theorem}\label{elt}
Let $B$ be a flat strip, then any geodesic in $B$ is closed and induce the same homotopy type. Moreover there are only finitely many homotopy classes of such geodesics. 
\end{theorem}
In fact, we conjecture that much more is true: there is a constant $C$ that depends only on the genus of $S$, the area of $S$, the minimum length of the geodesics in $S$ and the minimum of the curvature, such that the number of homotopy classes of such geodesics is bounded by $C$.

Finally, we warn the reader that we will work mostly in the universal covering, i.e. the disc, so, many times, we will call in the same way objects in the surface and in the universal covering. Any geodesic is parametrized by unit speed. For a geodesic we mean the whole geodesic, if we want to use a piece of geodesic we call it a geodesic segment.   
\end{section}

\begin{section}{\bf Preliminaries}\label{pre}
In this is section we will introduce the terminology and the definitions used in the rest of the paper. Let $S$ be a surface of genus bigger than one and let $p:\D\to S$ be the covering projection. Let us call $\Gamma$ the group of deck transformations so that $S=\D/\Gamma$. We may suppose taking an appropriated covering projection that $\Gamma$ is a subgroup of Moebius transformations of the unit disc. Let $g$ be a $\Cd^r$, $r\geq 2$, metric in $S$ and call also $g$ a lift of the metric to the disc.

Now we will list some known properties of such a subgroup $\Gamma$ and of nonnpositive curvature surfaces. For example, if $N_1,N_2\in\Gamma$ commute, i.e. $N_1\circ N_2=N_2\circ N_1$, then there is $M\in\Gamma$ and $k_1,k_2\in\Z$ such that $N_i=M^{k_i}$ for $i=1,2$. Also for any $M\in\Gamma$ there is a geodesic $\gamma$, associated to a closed geodesic of $S$, and a constant $c$ such that $M(\gamma(t))=\gamma(t+c)$ for any $t\in\R$. Conversely, given a geodesic $\gamma$ in $\D$ projecting to a closed geodesic of $S$ there is $M\in\Gamma$ and a constant $c$ such that $M(\gamma(t))=\gamma(t+c)$ for any $t\in\R$. So we get that $\gamma$ projects to a closed geodesic if and only if there is $M\in\Gamma$ and a constant $c$ such that $M(\gamma(t))=\gamma(t+c)$ for any $t\in\R$. On the other hand, every geodesic $\gamma$ has two endpoints in $\partial\D=S^1$, 
\begin{eqnarray*}
\gamma^+&=&\lim_{t\to +\infty}\gamma(t)\\
\gamma^-&=&\lim_{t\to -\infty}\gamma(t)
\end{eqnarray*}
with $\gamma^+\neq\gamma^-$ which uniquely determines $\gamma$ up to flat strips, i.e. if $\gamma_1^{\pm}=\gamma_2^{\pm}=\gamma^{\pm}$ for two geodesics $\gamma_1,\gamma_2$ then, as it is not hard to see they bound a region in $\D$. This region, which we shall call $B(\gamma_1,\gamma_2)=B$, satisfies $\K|_B\equiv 0$ and for any geodesic $\eta$ such that $\eta^-=\gamma^-$ we have that $\eta^+=\gamma^+$ and conversely. Moreover, each geodesic $\eta\subset B$ has its endpoints $\gamma^{\pm}$. In fact there is and isometry $h:L\to\overline{B}$ where 
$$
L=\{z\in\C\; \mbox{s.t.}\; 0\leq\re(z)\leq \alpha\}
$$ 
and $\alpha=d(\gamma_1,\gamma_2)$ sending $\re(z)=0$ to $\gamma_1$ and $\re(z)=\alpha$ to $\gamma_2$. So, using this isometry, we get that the area of $B$ equals $\infty$ and also equals $\infty$ any end of $B$, i.e. $h(\Im(z)>R)$ or $h(\Im(z)<R)$ any $R\in\R$. Moreover, given a flat strip $B$, if there is a $M\in\Gamma$ such that $M(B)\cap B\neq\emptyset$, then either $M(B)=B$ and moreover, $M(\gamma)=\gamma$ for any geodesic $\gamma\subset B$, or $M\gamma$ and $\gamma$ intersect at only one point and the angle they make at this point is not zero and do not depend on the geodesic $\gamma\subset B$. This fact is a concequence of the following more general fact and the flatness of the strip. The metric $g$ has no conjugated points, i.e. given two points $z,w\in\D$ there is one and only one geodesic $\gamma$ passing throw $z,w$. Moreover given $z\in S^1$ and $w\in\D$ there is one and only one geodesic $\gamma$ passing throw $w$ and beginning or ending at $z$ (beginning or ending depends on which direction we walk the geodesic, but the geodesic is unique). Given $M\in\Gamma$, we have that $M$ acts on $S^1$. $M$ has exactly two fixed points $z^+,z^-\in S^1$ which coinside with the endpoints of a geodesic $\gamma$ such that $M\gamma=\gamma$. Moreover, $M^n(z)\to z^+$ as $n\to +\infty$ for any $z\in\overline{\D}\minus\{z^-\}$ and $M^n(z)\to z^-$ as $n\to -\infty$ for any $z\in\overline{\D}\minus\{z^+\}$ where $\D=\D\cup\partial\D$.

\end{section}

\begin{section}{\bf Conformal metric}
In this section we shall prove that we may suppose that the metric induced in the disc is conformal, i.e. a scalar multiple of the euclidean metric. 
Let $S$ be a surface of genus bigger than one and let $p:\D\to S$ be the covering projection. Let us call $\Gamma$ the group of deck transformations so that $S=\D/\Gamma$. We may suppose taking an appropriated covering projection that $\Gamma$ is a subgroup of Moebius transformations of the unit disc. Let $g$ be a $\Cd^r$, $r\geq 2$, metric in $S$ and call also $g$ the lift of the metric to the disc. Given a map $\phi:\D\to\R$ and a subgroup of Moebius transformations of the unit disc $\Gamma'$ we say that $\phi$ is $\Gamma'-$equivariant if $\phi(M(z))|M'(z)|=\phi(z)$ for any $M\in\Gamma'$ and $z\in\D$. Let us denote with $<v,w>$ the standard inner product in $\R^2$. So, we have the following proposition:
\begin{prop}\label{conf}
There exists a subgroup $\Gamma'$ of Moebius transformations of the unit disc isomorphic to $\Gamma$ and a $\Cd^r$ positive $\Gamma'-$equivariant map $\rho:\D\to\R$ such that defining the metric $\hat{g}_x(v,w)=\rho^2(x)<v,w>$, there exists a $\Cd^{r+1}$ isometry $h:(\D,g)\to(\D,\hat{g})$ satisfying $\Gamma'=\{h\circ M\circ h^{-1}\;,\;M\in\Gamma\}$
\end{prop}
So we get that $S=\D/\Gamma$ with the metric $g$ is isometric with $\D/\Gamma'$ with the metric induced by $\hat{g}$. And so we will suppose that the metric $g$ is already given by the scalar $\rho^2$ times the euclidean metric. 

The idea of the proof of proposition \ref{conf} is that the metric $g$ defines a field of ellipses that should be sent to circles by the isometry, and this is done using Beltrami coefficients and the Measurable Riemann mapping theorem. So, let us begin with some definitions and properties (see \cite{lv} for the details). Given a linear map $L:\R^2\to\R^2$ there exists unique $a,b\in\C$ such that $L(v)=av+b\bar{v}$ for any $v\in\C$ where we denote with $\bar{v}$ the complex conjugate of $v$. We call $\mu_L=\frac{b}{a}$ the complex dilation of $L$. Now, for a differentiable map $h:U\to\C$ we have that $D_zh(v)=a(z)v+b(z)\bar{v}$ and we call $\partial h(z)=a(z)$ and $\bar{\partial} h(z)=b(z)$. When $h$ is a orientation preserving diffeomorphism, we define the Beltrami coefficient of $h$, or complex dilation of $Dh$ at $z$, by
$$
\mu_h(z)=\frac{\bar{\partial} h(z)}{\partial h(z)}
$$
As $h$ is an orientation preserving homeomorphism we get that $|\mu_h(z)|<1$ for any $z\in U$. $\mu_h$ essentially defines a field of ellipses, i.e. if we take the field of ellipses that are sent via $Dh$ to circles, then each ellipse is defined by the following: the ellipse at $z$ has eccentricity $\frac{1+|\mu_h(z)|}{1-|\mu_h(z)|}$ and its mayor axis is $\sqrt{\mu_h(z)}$ whenever $\mu_h(z)\neq 0$. If $\mu_h(z)=0$ then $D_zh$ sends circles into circles, i.e. it is conformal. Also we have a converse to this construction. 

\begin{theorem}{\bf Measurable Riemann mapping theorem} 
Given $\mu:U\to\C$ measurable such that $\|\mu\|_{\infty}<1$ where $\|\cdot\|_{\infty}$ is the $L^{\infty}$ norm, there is a quasiconformal homeomorphism $h:U\to\C$ such that $\mu=\mu_h$ a.e. and $h$ is unique modulo composition with a holomorphic map, i.e. if there is $h'$ such that $\mu_{h'}=\mu$ a.e. then there is an holomorphic $f:h(U)\to\C$ such that $h'=f\circ h$. Moreover if $\mu=0$ a.e., $h$ is holomorphic and if $\mu$ is $\Cd^r$ then $h$ is a $\Cd^{r+1}$ diffeomorphism.   
\end{theorem}

Now, we say that a Beltrami coefficient is $\Gamma-$equivariant for some subgroup $\Gamma$ of Moebius transformations if 
$$
M_{*}\mu=\mu\circ M\frac{\overline{M'}}{M'}=\mu
$$
for any $M\in\Gamma$. It is clear from the definition that if $h$ is a solution to $\mu_h=\mu$ and $\mu$ is $\Gamma-$equivariant then $h\circ M$ is also a solution to the equation for any $M\in\Gamma$ and so, by uniqueness, there are $f_M$ holomorphic such that $h\circ M=f_M\circ h$. So we are ready to prove the proposition.

\begin{proof}{\bf of Proposition \ref{conf}}
Given the metric $g$, there are positive definite symmetric matrices $A(z)$ for $z\in\D$ such that $g_z(v,w)=<A(z)v,w>$. As $g$ comes from the surface we get that 
$$
\overline {M'(z)}A(M(z))M'(z)=A(z)
$$ 
for any $M\in\Gamma$. We seek for an $h:\D\to\D$ and a map $\rho:\D\to\R$ such that 
\begin{eqnarray}\label{eq1}
<A(z)v,w>=\rho^2(h(z))<(D_zh)^tD_zhv,w>
\end{eqnarray}
for any $z\in\D$ and $v,w\in\R^2$, where $(D_zh)^t$ is the adjoint of $D_zh$. But equation (\ref{eq1}) means that $A(z)=\rho^2(h(z))(D_zh)^tD_zh$. Now if we take $\mu_A(z)$ the complex dilation of $A$ at $z$ then an easy calculation gives that we need  
$$
\mu_h(z)=\frac{\mu_A(z)}{|\mu_A(z)|^2}\Bigl(1-\sqrt{1-|\mu_A(z)|^2}\Bigr)  
$$
whenever $\mu_A(z)\neq 0$. The equivariance of the metric $g$ and so the one of $A$ gives us that $M_*\mu_A=\mu_A$ for any $M\in\Gamma$ and so if we define
$$
\mu(z)=\frac{\mu_A(z)}{|\mu_A(z)|^2}\Bigl(1-\sqrt{1-|\mu_A(z)|^2}\Bigr)  
$$
whenever $\mu_A(z)\neq 0$ and $\mu(z)=0$ if $\mu_A(z)=0$, we get that $M_*\mu=\mu$ for any $M\in\Gamma$ and so $\mu$ is $\Gamma-$equivariant. Now, as $\mu_A$ is $\Gamma-$equivariant and $\D/\Gamma=S$ is compact, taking a compact fundamental domain and noting that $|M_*\mu|=|\mu\circ M|$, it is not hard to see that $\|\mu_A\|_{\infty}=n<1$, it is essentially because the eccentricities of the ellipses are uniformly bounded by the compactness of the surface. So we get that 
$$
\|\mu\|_{\infty}<\frac{1-\sqrt{1-n^2}}{n}<n<1
$$
So, using the measurable Riemann mapping theorem, we get that there is a quasiconformal homeomorphism $h:\D\to\C$ such that $\mu_h=\mu$. Moreover, $h$ is $\Cd^{r+1}$ if $\mu$ is $\Cd^r$. As $h(\D)$ must be open and simply connected and $h(\D)\neq\C$, after composing with a biholomorphism we may suppose that $h(\D)=\D$. $\mu$ being $\Gamma-$equivariant implies that $h\circ M:\D\to\D$ is also a solution to the Beltrami equation for any $M\in\Gamma$, an so defining 
$$
\Gamma'=\{h\circ M\circ h^{-1} \;\mbox{s.t.}\;M\in\Gamma\}
$$
we get that $\Gamma'$ is a subgroup of Moebius transformations of the unit disc which is isomorphic to $\Gamma$. By construction, we get that the complex dilation of $A(z)$ is equal to the one of $(D_zh)^tD_zh$. Both being positive definite selfadjoint matrices, it is not hard to see that there exists a positive number $\rho^2(h(z))$ such that $A(z)=\rho^2(h(z))(D_zh)^tD_zh$. The differentiability of $\rho$ follows from the differentiability of $A$ (i.e. the metric $g$) and the one of $h$. also it follows the equivariance of $\rho$. As we saw, for this $\rho$, equation (\ref{eq1}) is satisfied, and so we proved the proposition.
\end{proof}

So, from now on we suppose that $S=\D/\Gamma$ for some group $\Gamma$ of Moebius transformations and that the metric is given by a positive map $\rho:\D\to\R$ satisfying $\rho(M(z))|M'(z)|=\rho(z)$. We call $p:\D\to S$ the covering projection. We point out that for such metric the curvature is given by
$$
\K(z)=-\frac{\Delta\log\rho(z)}{\rho^2(z)}
$$
where $\Delta\log\rho$ is the laplacian of $\log\rho$.
\end{section}
\begin{section}{\bf Isometries of the plane}
We shall use some results on subgroups of orientation preserving isometries of the euclidean plane. So we will give some properties of them. Recall that the subgroup of orientation preserving euclidean isometries of the plane is
$$
Iso(\C)=\{\lambda z+a\;\mbox{s.t.}\;|\lambda|=1\;\mbox{and}\;a\in\C\}
$$ 
When we say that a subgroup $G$ is isomorphic to a subgoup $G'$ we mean that there is $L\in Iso(\C)$ such that 
$$
LGL^{-1}=\{L\circ R\circ L^{-1}\;\mbox{s.t.}\;R\in G\}=G'
$$
  
Let $G\subset Iso(\C)$ be a subgroup of isometries of the plane. We say that a set $A$ is $G-$invariant if $R(A)=A$ for any $R\in G$. We say that $G$ is minimal if the only open invariant sets are the whole plane and the emptyset. 

We call $G$ a subgroup of rotations if it is isomorphic to a subgroup of 
$$
Rot(\C)=\{\lambda z\;\mbox{s.t.}\;|\lambda|=1\}
$$  
Recall that a subgroup of rotations acts minimally on circles centered at the origin of the rotation (i.e. the only open invariant subsets of a circle centered at the origin of the rotation are the emptyset and the circle itself) or else it is finite. 

Let us call
$$
Trans(\C)=\{z+a\;\mbox{s.t.}\;a\in\C\}
$$
We have that $Trans(\C)\simeq\C$. So we have that any subgroup $G$ of $Trans(\C)$ is discrete or not. If it is discrete, it may have one generator and there exist $\alpha\in\R$ such that $G$ is isomorphic to 
$$
\Z(\alpha)=\{z+n\alpha\;\mbox{s.t.}\;n\in\Z\}
$$ 
or two generators and hence there exist $\alpha\in\R$ and $a\in\Hy$ where $\Hy=\{z\in\C\;\mbox{s.t.}\;\Im(z)>0\}$ such that $G$ is isomorphic to 
$$
\Z^2(\alpha,a)=\{z+n\alpha+ma\;\mbox{s.t.}\;n,m\in\Z\}
$$
If it is not discrete, then either it is minimal or it acts minimaly on lines parallel to an axis of the group, i.e. there exists $0\neq a\in\C$ and $\alpha\in\R\minus\Q$ such that 
$$
G\supset\{z+(n+m\alpha)a\;\mbox{s.t.}\;m,n\in\Z\}
$$ 
and also 
$$
G\subset\{z+ra\;\mbox{s.t.}\;r\in\R\}
$$  

For $\alpha\in\R$ let us define the following subgroups of $Iso(\C)$ 
\begin{eqnarray*}
\Lambda_0(\alpha)&=&\{e^{\frac{\pi}{3}ki}z+(n+e^{\frac{\pi}{3}i}m)\alpha\;\;\mbox{s.t.}\;k,n,m\in\Z\}\\
\Lambda_1(\alpha)&=&\{e^{\frac{2\pi}{3}ki}z+(n+e^{\frac{\pi}{3}i}m)\alpha\;\;\mbox{s.t.}\;k,n,m\in\Z\}\\
\Z^2_i(\alpha)&=&\{i^kz+(n+im)\alpha\;\;\mbox{s.t.}\;k,n,m\in\Z\}\\
\Z_-(\alpha)&=&\{(-1)^kz+n\alpha\;\;\mbox{s.t.}\;k,n\in\Z\}
\end{eqnarray*}
Finally for $\alpha\in\R$ and $a\in\Hy$ let us put
$$
\Z^2_-(\alpha,a)=\{(-1)^kz+n\alpha+ma\;\;\mbox{s.t.}\;k,n,m\in\Z\}
$$
Let us call these subgroups exceptional. We claim that all the list above exhaust the subgroups of orientation preserving isometries of the plane up to isomorphisms.
\begin{prop}\label{clas}
Let $G\subset Iso(\C)$ be a subgroups of orientation preserving isometries of the plane. Then one of the following hold:
\begin{enumerate}
\item $G$ is minimal.
\item $G$ is a soubgroup of rotations acting minimally on cicles centered at the origin of the rotation. 
\item $G$ is a finite soubgroup of rotations and so there exists $k\in\Z$ such that $R^k=id$ for any $R\in G$
\item $G$ acts minimally on lines parallel to an axis of the group.
\item $G$ is isomorphic to one of the form $\Z(\alpha)$.
\item $G$ is isomorphic to one of the form $\Z^2(\alpha,a)$.
\item $G$ is isomorphic to an exceptional subgroup.
\end{enumerate}
\end{prop}
For the proof of the proposition let us define for a subgroup $G$, a homorphism $d:G\to S^1$ by $d(T)=T'(0)$. Notice that $\ker d=G\cap Trans(\C)$.
\begin{proof}
If $\ker d=\{id\}$ Then $G$ is a soubgroup of rotations and hence it is clear that 2) or 3) must hold. So let us suppose $\ker d\neq\{id\}$. Now 
if $\ker d$ is not discrete we get case 1) or case 4). So let us suppose that $\ker d$ is discrete. Suppose first that $\ker d$ has one generator and call it $a$, so that $\ker d=\{z+na\;\mbox{s.t.}\;n\in\Z\}$. We claim that $\im d\subset\{\pm id\}$ and so we get case 5) or 7). To proof our claim suppose that there is $1\neq\lambda\in\im d$ conjugating with a translation we may suppose that $R(z)=\lambda z$ and $T(z)=z+a$ belong to $G$. So $R\circ T\circ R^{-1}=z+\lambda a$ also belong to $\ker d$ and so $\lambda a=na$ for some $n\in\Z$ hence $\lambda=-1$. In case $\ker d$ has two generators, we claim that for any $\lambda\in\im d$ we have that $\lambda^6=1$ or $\lambda^4=1$. Let us call $a$ and $b$ the two generators so that  $\ker d=\{z+na+mb\;\mbox{s.t.}\;n,m\in\Z\}$ and take $1\neq\lambda\in\im d$. Conjugating with a translation we may suppose that $R(z)=\lambda z$ belong to $G$ and $T_a=z+a$ and $T_b=z+b$ still are the generators of $\ker d$. So $R\circ T_{\sigma}\circ R^{-1}=z+\lambda\sigma$ for $\sigma=a,b$ also belong to $\ker d$ and so there are $r_{\sigma},s_{\sigma}\in\Z$, $\sigma=a,b$ such that
\begin{eqnarray*}
\lambda a&=&r_aa+s_ab\\
\lambda b&=&r_ba+s_bb
\end{eqnarray*}
This implies that $\lambda$ must be a root of
$$
\lambda^2-(r_a+s_b)\lambda+r_as_b-r_bs_a=0
$$
which is a polinomial with integer coefficients of degree two. If $\lambda =-1$ we have nothing to do so we may suppose $\Im \lambda\neq 0$ and so as $\bar{\lambda}$ must also be a root of the polinomial and $|\lambda|=1$ we get that $r_as_b-r_bs_a=1$. Also we get that $|r_a+s_b|\leq 2$ so the polinomial transform into
$$
\lambda^2-n\lambda+1=0
$$
with $n=-2,-1,0,1,2$. Which implies $\lambda^6=1$ or $\lambda^4=1$. Now as $\im d$ is a subgroup of $S^1$ we get that 
$$
\im d\subset G_6=\{\lambda\;\mbox{s.t}\;\lambda^6=1\}
$$
or
$$
\im d\subset G_4=\{\lambda\;\mbox{s.t}\;\lambda^4=1\}
$$
From this it is not hard to see that $G$ must be exceptional or we are in case 6).
\end{proof}

In cases 5,6,7) we have that the action acts discontinuously, maybe with fixed points. So we can make the quotient of $\C\minus\mbox{Fix}(G)$ by the action. In case 5) we get a cilinder, in case 6) a torus and in case 7) we get the sphere with some points removed.

\end{section}
\begin{section}{\bf Flat strips}\label{flat}
Define
$$
\hat W=\{x\in\D\;\;\mbox{s.t. there is an open flat strip}\;B\;\mbox{with}\;x\in B\} 
$$
Clearly $\Gamma(\hat W)=\hat W$, i.e. $M(\hat W)=\hat W$ for any $M\in\Gamma$ and $\hat W$ is open. Call $W$ a connected component of $\hat W$ and call $\Gamma_0$ the largest subgroup of $\Gamma$ satisfying $\Gamma_0(W)=W$ i.e.
$$
\Gamma_0=\{M\in\Gamma\;\mbox{s.t.}\;M(W)=W\}
$$
Notice that $M(W)\cap W\neq\emptyset$ implies $M(W)=W$. It is not hard to see that any bounded connected component of the complement of $W$ is a convex set with boundary a geodesic poligonal consisting of fintely many geodesic segments. Moreover, the boundary of the unbounded components are piecewise geodesic curves which descends to a simple closed curve with finetely many geodesic segments and homotopic to a simple closed geodesic. 

In this section we will prove some properties of $W/\Gamma_0$ in order to get the proof of the theorem in the next section. So let us start.

Take $g:\Hy\to W$ an holomorphic covering map and define $\hat{\rho}$ the pullback of the conformal metric, i.e. $\hat{\rho}(z)=|g'(z)|\rho(g(z))$. As 
$$
\K(z)=-\frac{\Delta\log\rho(z)}{\rho^2}=0
$$
on $W$ we get that $\log\rho$ is harmonic on $W$. Also, as
$$
\hat{\K}(z)=-\frac{\Delta\log\hat{\rho}(z)}{\hat{\rho}^2}=\K(g(z))=0
$$
we get that $\log\hat{\rho}$ is harmonic on $\Hy$. As $\Hy$ is simply connected we have $h:\Hy\to\C$ holomorphic such that $\log |h'(z)|=\log\hat{\rho}(z)$ or $|h'(z)|=\hat{\rho}(z)$. So we have that $h$ is a local isometry of the conformal metric $\hat{\rho}$ on $\Hy$ to the euclidean metric on $\C$. Call $E=h(\Hy)$. As $h$ is holomorphic $E$ is open. Notice that $h'(z)\neq 0$.

Now call $\Gamma_1$ the pullback of $\Gamma_0$ by $g$, i.e.
$$
\Gamma_1=\{M\in\mbox{Moebius}(\Hy)\;\;\mbox{s.t. there is}\; N\in\Gamma_0\;\mbox{with}\;g\circ M=N\circ g\}
$$
and call $G$ the pushforward by $h$ of $\Gamma_1$, i.e
$$
G=\{R\in Iso(\C)\;\mbox{s.t.}\;h\circ M=R\circ h\;\mbox{for some}\;M\in\Gamma_1\}
$$
Now, for $G$ we have that Proposition \ref{clas} holds, so one of the seven properties must hold. We claim that only properties 5), 6) or 7) can hold. Let us see how this claim follows. 1) cannot hold because if $G$ is minimal, as $E$ is open and invariant by $G$, $E$ must be equal to $\C$ and there is no holomorphic map from $\Hy$ onto $\C$ with $h'\neq 0$. 

\begin{subsection}{\bf Case 2)}
Let us see that property 2) cannot hold. Suppose $G$ is a subgroup of rotations centered at the origin (just compose $h$ with a translation) acting minimally on circles centered at $0$. We claim that there is an $r\geq 0$ such that $E=\C\minus\overline{\D}_r$ where $\overline{\D}_r=\{z\in\Z\;\mbox{s.t.}\;|z|\leq r\}$. Take $r=\inf\{|z|\;\mbox{s.t.}\;z\in E\}$ then clearly $E\subset\C\minus\overline{\D}_r$. Now take $z$ with $|z|>r$. As $E$ is unbounded as it must contain a line (coming from a flat strip in $W$) and connected, there must be $w\in E$ with $|w|=|z|$. Now, $E$ is open, so we have that there is a disk centered at $w$ of positive radius. The intersection of this disk with the circle centered at $0$ and radius $|z|$ is an open arc of the circle and hence as $G$ acts minimally on it and $E$ is $G-$invariant we get that the whole circle is in $E$ and hence $z\in E$.  Now $r$ must be positive because there is no holomorphic map from $\Hy\to\C\minus\{0\}$ without $0$ derivative. As $\hat h:\Hy\to\C\minus\overline{\D}_r$ defined by $\hat h(z)=re^{-iz}$ is a covering map we get that $h$ equals to $\hat h$ modulo composition with a Moebius transformation of $\Hy$. So changing $h$ and $g$ we may suppose that $h(z)=re^{-iz}$ and all the other properties still hold. But now, it is not hard to see that 
$$
\Gamma_1\subset\{z+\alpha\;\mbox{s.t.}\;\alpha\in\R\}
$$
Notice that as $g$ is a covering map and the definition of $\Gamma_1$ we have that $g(x)=g(y)$ implies that there is $N\in\Gamma_1$ such that $N(x)=y$. Moreover, $g(x)=Mg(y)$ for some $M\in\Gamma_0$ if and only if there is $N\in\Gamma_1$ such that $N(x)=y$. So we have that 
$$
\hat g:(\Hy/\Gamma_1,\hat \rho)\to (W/\Gamma_0,\rho)
$$
is a biyective isometry where $\hat g$ is the map induced by $g$. But then, it is not hard to see that $\Gamma_1$ acts discontinuously on $\Hy$ and also that we have
\begin{eqnarray}\label{area1}
\mbox{Area}(\Hy/\Gamma_1)=\mbox{Area}(W/\Gamma_0)\leq\mbox{Area}(S)
\end{eqnarray} 
As $\Gamma_1$ acts discontiuously we get that $\Gamma_1=\{z+n\beta\;\mbox{s.t.}\;n\in\Z\}$ for some $0<\beta\in\R$. As $\hat\rho(z)=|h'(z)|$ we have that $\hat\rho(z)=re^{\Im{z}}$. So as $0\leq\re(z)\leq\beta$ is a fundamental domain of $\Hy/\Gamma_1$ we get that it must have infinite area contradicting (\ref{area1}). So we finished with case 2).
\end{subsection}

\begin{subsection}{\bf Case 3)}
Let us see that property 3) cannot hold.  The idea is that we can take a finite covering of  $S$ so that $G$ becomes trivial, which implies that the corresponding $\hat\Gamma_0$ becomes trivial also wich contradicts the existence of a flat strip which has $\infty$ measure. Let us see how to do this. Suppose the rotations are centered at $0$ and call $k_0$ the order of $G$. Call $h_*:\Gamma_1\to G$ and $g_*:\Gamma_1\to\Gamma_0$ the corresponding actions of $h$ and $g$ on the groups i.e. $h_*(M)=R$ if $h\circ M=R\circ h$ and the same with $g$. We have that 
$$
\Gamma_1/\!\ker h_*\simeq G
$$
We claim that the area of $\Hy/\!\ker h_*$ is finite (the area induced by $\hat\rho$). In fact, we have a covering map
$$
\Hy/\!\ker h_*\to\Hy/\Gamma_1
$$
which is $k_0$ to $1$. Then we get that $\mbox{Area}(\Hy/\!\ker h_*)=k_0\mbox{Area}(\Hy/\Gamma_1)$. By equality (\ref{area1}) we get our claim. Now take a a flat strip $B\subset W$ and lift it to $\hat B\subset\Hy$ by $g$. We get that $\hat B$ is also a flat strip. Let us state the following usefull lemma:
\begin{lema}\label{int1}
For any flat strip $\hat B\subset \Hy$ there is $id\neq M\in\ker h_*$ such that $M(\hat B)\cap\hat B\neq\emptyset$.
\end{lema} 
\begin{proof}
We have that the area of $\hat B$ is $\infty$ and so it cannot decend injectibly to $\Hy/\!\ker h_*$ which has finite area, so there must be $id\neq M\in\ker h_*$ such that $M(\hat B)\cap \hat B\neq\emptyset$.
\end{proof}
Take a geodesic $\gamma\subset\hat B$ then we have that $M\gamma\cap\gamma\neq\emptyset$. So we have to possibilities, either \begin{enumerate}[i)]
\item $M\gamma\cap\gamma=\{x\}$ and they make an angle $\theta\neq 0$ at $x$ or 
\item $M\gamma=\gamma$ and there is $c>0$ such that $M(\gamma(t))=\gamma(t+c)$ for any $t\in\R$.
\end{enumerate}
Notice that $h(\gamma(t))=at+b$ and $h(M(\gamma(t)))=\hat at+\hat b$ are two lines $l_1,l_2$ respectively. Now, in case i) we get a contradiction because $l_1$ and $l_2$ must be two lines intersecting at $h(x)$ and making angle $\theta$ at $h(x)$ because $h$ is conformal, but on the other hand, as $M\in\ker h_*$ we get that $h\circ M=h$ and so $l_2=h(M\gamma)=h(\gamma)=l_1$. In case ii) we get a contradiction also because $h(\gamma(t+c))=h(M\gamma(t))=h(\gamma(t))$ for any $t\in\R$ and so, $at+b=a(t+c)+b$ for any $t\in\R$ which is clearly a contradiction whith $a,c\neq 0$. So we are done with case 3).
\end{subsection}

\begin{subsection}{\bf Case 4)}
Let us see that property 4) cannot hold. Suppose, composing $h$ with a rotation of the plane if necesary that 
$$
G\subset\{z+r\;\mbox{s.t.}\;r\in\R\})
$$  
So we get, as $E$ is open and $G$ acts minimally on lines paralell to the real axis, that 
$$
E=\{z\in\C\;\mbox{s.t.}\;\alpha<\Im z<\beta\}
$$
for some $\alpha$ and $\beta$ in $\R$, maybe $\pm \infty$. Notice that $E\neq\C$. Now, $E$ is simply connected, so $h$ is a biholomorphism. As $G$ does not act discontinously we get that $\Gamma_1$ does not act discontinuously on $\Hy$ contradicting what we saw in case 2). 
\end{subsection}

So we get that only cases 5), 6) or 7) can hold. 

\end{section}

\begin{section}{\bf Prove of Theorem \ref{elt}}
In this section we will prove the Theorem \ref{elt}, i.e. that any geodesic contained in a flat strip represent a closed geodesic (it is clear from the properties listed on the section of preliminaries that two such geodesics induce the same homotopy type) and that there are only finitely homotopy clases of them. We have defined $\hat W$ as the set of points belonging to an open flat strip. So we are going to prove that any geodesic in a connected component $W$ of $\hat W$ represents a closed geodesic and that there are only finetely many clases for each component. Once we get this, the whole result follows because the number of connected components cannot exced $3g-3$ (g=genus of $S$) because each one contains a closed geodesic. 

Fix a connected component $W$. We use the same terminology as in the preceding section. So we have to deal with cases 5), 6) and 7) of Proposition \ref{clas}. 
\begin{subsection}{\bf Case 5)}\label{cili}
In this case we have that $G$ is isomorphic to 
$$
\Z(\alpha)=\{z+n\alpha\;\;\mbox{s.t.}\;n\in\Z\}
$$
for some $\alpha\in\R$. So suppose $G=\Z(\alpha)$. We get that $E/G$ is a subset of the cilinder $\C/\Z(\alpha)$. Let us prove first that any geodesic $\gamma\subset B\subset W$ in a flat strip represents a closed geodesic. To do this let us state an analogous of lemma \ref {int1}:
\begin{lema}
For any flat strip $B\subset W$ there is $id\neq M\in\Gamma_0$ such that $M(B)\cap B\neq\emptyset$.
\end{lema} 
\begin{proof}
We have that the area of $B$ is $\infty$ and so it cannot decend injectibly to the surface, so there must be $M\in\Gamma$ such that $M(B)\cap B\neq\emptyset$, but by the definition of $\Gamma_0$ this implies that $M\in\Gamma_0$.
\end{proof}
So, as $\gamma\subset B$ for some flat strip $B$, we get that either $M\gamma=\gamma$ and hence $\gamma$ represents a closed geodesic or $M\gamma\cap\gamma=\{x_0\}$ and the angle between $M\gamma$ and $\gamma$ at $x_0$ is $\theta\neq 0$. Lift $\gamma$, $M\gamma$ and $x_0$ to $\gamma_i$, $i=1,2$ and $\{x\}=\gamma_1\cap\gamma_2$ by $g$. As $g(\gamma_2)=Mg(\gamma_1)=g(N_0\gamma_1)$ for some $N_0$ such that $M\circ g=g\circ N_0$ we get that $N\gamma_1=\gamma_2$ for some $N\in\Gamma_1$. The angle between $\gamma_1$ and $N\gamma_1$ at $x$ is again $\theta$ because $g$ is conformal. Sending $\gamma_1$ to $E$ by $h$ we get a line $l=h(\gamma_1)$, $z=h(x)$ and $n\in\Z$ such that $l+n\alpha\cap l=\{z\}$. But $l$ and $l+n\alpha$ are paralell lines, so this is a contradiction. So we get that any geodesic in $W$ represents a closed geodesic. Moreover, the same reasoning, using that $G$ has only one generator, prove that if $\gamma_1$ and $\gamma_2$ are two geodesics in $W$ then they cannot intersect. This implies that $W$ is convex, in fact, take two points in $W$, then, if they do  not belong to a same flat strip, it is not hard to build a curve $\phi$ consisting of a finite number of legs, each leg belonging to a geodesic contained in $W$. But as we saw, two distinct geodesics cannot intersect, so the curve consist of only one geodesic arc, and so $W$ is convex. So $W$ is simply connected and hence $g$ is a biholomorphism. Moreover, $f=h\circ g^{-1}$ must be also a biholomorphism by the convexity of $W$. So we get that $E$ is a strip between to lines and hence $W$ consist of only one flat strip and so $W/\Gamma_0$ is isometric with a flat cilinder of finite area. This implies that there is only one homotopy type of geodesic in $W/\Gamma_0$.
\end{subsection}

\begin{subsection}{\bf Case 6)}\label{tor}
In this case we have that $G$ is isomorphic to 
$$
\Z^2(\alpha,a)=\{z+n\alpha+ma\;\;\mbox{s.t.}\;n,m\in\Z\}
$$
for some $\alpha\in\R$ and $a\in\Hy$. So suppose $G=\Z^2(\alpha,a)$.
The same prove of case 5) gives that any geodesic in $W$ represents a closed geodesic. Notice that $E\neq\C$. Moreover it is not hard to see that the complement of $E$, $E^c$ has nonempty interior. So we have that $E/G\subset\C/\Z^2(\alpha,a)$ and the inclusion is strict. We make the following two claims:
\begin{claim}\label {3g}
For any direction $\sigma\in S^1$ there are at most $3g-3$ homotopy classes of geodesics in $W$ satisfying the following: if we take $\gamma\subset W$ one of this geodesics, and lift it to $\eta\subset\Hy$ by $g$ so that $g\circ\eta=\gamma$ then $h\circ\eta$ is a line with direction $\sigma$. 
\end{claim}

\begin{claim}\label{fini}
The number of directions $\sigma\in S^1$ for which there is a geodesic $\gamma\subset W$ such that lifting it to $\eta\subset\Hy$ by $g$ so that $g\circ\eta=\gamma$ then $h\circ\eta$ is a line with direction $\sigma$, is finite.
\end{claim}
The Theorem is straightforward in this case using the two claims.
\begin{proof}{\bf of Claim \ref {3g}}
Take a direction $\sigma$ and $\gamma_i\subset W$, $i=1,2$ two geodesics such that lifting them to $\eta_i\subset\Hy$ by $g$ so that $g\circ\eta_i=\gamma_i$ then $h\circ\eta_i$ is a line with direction $\sigma$. We claim that either $\gamma_1\cap\gamma_2=\emptyset$ or $\gamma_1=\gamma_2$. Suppose that $\gamma_1\cap\gamma_2\neq\emptyset$ and that they make an angle $\theta\neq 0$ at the point of itersection. We can lift them to $\hat{\eta}_1,\eta_2\subset\Hy$ in order to get $\hat{\eta}_1\cap\eta_2\neq\emptyset$. Moreover, as $\hat{\eta}_1$ and $\eta_1$ are two lift of the same geodesic we get that there is $N\in\Gamma_1$ such that $N\hat{\eta}_1=\eta_1$. As $h\circ N=T\circ h$ for some translation $T$ we get that $h\circ\hat{\eta}_1$ and $h\circ\eta_1$ have the same direction. So let us assume $\eta_1\cap\eta_2\neq\emptyset$. As $g$ is conformal we get that $\eta_1$ and $\eta_2$ have also angle $\theta$ at the point of intersection. But now, as $h$ is also conformal, $h\circ\eta_1$ and $h\circ\eta_2$ have angle $\theta$ at the point of intersection contradicting they have the same direction. So the claim follows from the fact there are at most $3g-3$ pairwise disjoint non homotopic closed geodesics in a surface of genus $g$. 
\end{proof}
\begin{proof}{\bf of Claim \ref {fini}}
Given a direction $\sigma$ in the hypothesis of the claim we get a line in $E$ with this direction, and hence a line in the open set $E/G\subset\C/\Z^2(\alpha,a)$ which has complement with nonempty interior, so we have that $\sigma$ must be rational, i.e. 
$$
\sigma=\frac{p\alpha+qa}{|p\alpha+qa|}
$$
where $p,q\in\Z$, $p,q$ coprime. Take a disc $D$ of radius $r$ in the complement of $E$. Then as there is a line whith direction $\sigma$ in $E$ we get that if we take all the lines whith direction $\sigma$ throu $D$ and all its translates $D+n\alpha+ma$, $n,m\in\Z$, they cannot fill all the plane. An easy calculation gives that this implies that 
$$
|p\alpha+qa|\leq\frac{\mbox{Area}(\C/\Z^2(\alpha,a))}{r}
$$
where $\mbox{Area}(\C/\Z^2(\alpha,a))=\alpha\Im a$.
\end{proof}
So we get the Theorem in this case. Notice that we proved also that the geodesics are simple closed geodesics.
\end{subsection}

\begin{subsection}{\bf Case 7)}
In this case we have that $G$ is exceptional, so it contain a subgroup of rotations of order less than $6$. The case $G=\Z_-(\alpha)$ follows exactly as Case 5). The case $G=\Z^2_-(\alpha,a)$ follows exactly as Case 6). So we deal with the other cases.

Recall the definition of the homorphism $d:G\to S^1$ by $d(T)=T'(0)$. Notice that $\ker d=G\cap Trans(\C)$, call $T(G)=\ker d$. Recall also the definition of $g_*:\Gamma_1\to\Gamma_0$ and $h_*:\Gamma_1\to G$. Call $T(\Gamma_1)=h_*^{-1}(T(G))$. We have that $G/T(G)=\im d$ which has order less than $6$ and $G/T(G)\simeq\Gamma_1/T(\Gamma_1)$. So we get a covering
$$
\Hy/T(\Gamma_1)\to\Hy/\Gamma_1
$$
which is at most $6$ to $1$. From this it follows that 
$$
\mbox{Area}(\Hy/T(\Gamma_1))\leq 6\mbox{Area}(\Hy/\Gamma_1)
$$
Now using the group $T(\Gamma_1)$ it is not hard to mimic the proof of Case 5) that any geodesic in a flat strip in $W$ is a closed geodesic. Moreover,  
Claim \ref{fini} of Case 6) applies as well. So we only have to prove that for any direction it correspond a finite number of homotopy classes of geodesics in $W/\Gamma_0$. The same proof of Claim \ref{3g} in Case 6) gives us the following claim
\begin{claim}
Given a direction $\sigma\in S^1$ and $\eta_1,\eta_2\subset\Hy$ such that $h\circ\eta_i$, $i=1,2$ have direction $\sigma$, we have that either $\eta_1\cap\eta_2=\emptyset$ or $\eta_1=\eta_2$. 
\end{claim}
So as in Case 6) we get that the geodesics in $\Hy/T(\Gamma_1)$ corresponding to a fix direction $\sigma$ are pairwise disjoint. Now we have the covering 
$$
\Hy/T(\Gamma_1)\to\Hy/\Gamma_1\simeq W/\Gamma_0
$$
which is at most $6$ to $1$. So as in $W/\Gamma_0$ there are at most $3g-3$ pairwise disjoint non-homotopic geodesics, it is not hard to see that the number of homotopy classes of geodesics in $\Hy/T(\Gamma_1)$ corresponding to $\sigma$ is finite for any $\sigma$. So there are at most finitely many homotopy classes of geodesics in $\Hy/T(\Gamma_1)$ and again using the covering, we get the theorem.

\end{subsection}

\end{section}

\bibliographystyle{alpha}

\end{document}